\documentclass[]{article}
\usepackage{amsmath}
\usepackage{amssymb}
\usepackage[margin=2cm]{geometry}
\usepackage{hyperref}
\usepackage{pgf,pgfarrows,pgfnodes,pgfautomata,pgfheaps,pgfshade}
\usepackage{graphicx,amsfonts}
\usepackage{amsmath,amssymb,setspace}
\usepackage{mathrsfs}
\usepackage{colortbl}
\usepackage{tikz}
\usepackage{framed}
\usepackage{subfig}
\usepackage[section]{placeins}
\usepackage{float}
\usepackage{enumerate}

\hypersetup{
    colorlinks,
    urlcolor=blue
}

\numberwithin{equation}{section}
\newtheorem{thm}{Theorem}[section]

\newtheorem{definition}[thm]{Definition}

\title{Relatively Coherent Sets as a Hierarchical Partition Method}

\author{Tian Ma and Erik M. Bollt}

\date{\today}

\begin{document}

\maketitle

%%%%%%%%%%%%%%%%%%%%%%%%%%%%%%%%%%%%%%%%%%%%%%%%%%%%%%%%%%%%%%%%%%%%%%%%%%%%%%%%%%
%                                                                                                                                      Abstract
%%%%%%%%%%%%%%%%%%%%%%%%%%%%%%%%%%%%%%%%%%%%%%%%%%%%%%%%%%%%%%%%%%%%%%%%%%%%%%%%%%

\begin{abstract}
\setcounter{section}{0}

Finite time coherent sets \cite{FSM} have recently been defined by a measure based objective function describing the degree that sets hold together, along with a Frobenius-Perron transfer operator method  to produce optimally coherent sets. 
Here we present an extension to generalize the concept to hierarchically define {\t relatively coherent sets} based on adjusting the finite time coherent sets to use relative measure restricted to sets which are developed iteratively and hierarchically in a tree of partitions.
Several examples  help clarify the meaning and expectation of the techniques, as they are the nonautonomous double gyre, the standard map, an idealized stratospheric flow, and empirical data from the Mexico Gulf during the 2010 oil spill.
Also for sake of analysis of computational complexity, we include an appendix concerning the computational complexity of developing the Ulam-Galerkin matrix estimates of the Frobenius-Perron operator centrally used here.

\vspace{8pt}
\noindent
{\small \bf Keywords: Coherent Pairs, Relative Measures, Relatively Coherent Structures, Frobenius-Perron Operator, Subdivision Flow Chart.} 
\vskip .2 cm
\noindent

\end{abstract}
%\FloatBarrier

%%%%%%%%%%%%%%%%%%%%%%%%%%%%%%%%%%%%%%%%%%%%%%%%%%%%%%%%%%%%%%%%%%%%%%%%%%%%%%%%%
%                                                                                                                                   1 Introduction
%%%%%%%%%%%%%%%%%%%%%%%%%%%%%%%%%%%%%%%%%%%%%%%%%%%%%%%%%%%%%%%%%%%%%%%%%%%%%%%%%
\section{Introduction}

Central to understanding mixing and transport mechanisms is the related question of defining partitions relative to which the transport can be discussed.  To this end, the concept of almost invariant sets defined for autonomous systems \cite{DJ, BN, FK} and coherent sets for non autonomous systems \cite{FSM} are central since transport may be defined as the measure of the set that leaves the partition element corresponding to the almost invariant set (or finite time coherent set) in a given time epoch. 
See also \cite{FK}.
Transfer operator methods are proven to be computationally effective for use in identifying almost invariant sets and finite time coherent sets. 
See \cite{BLSB, FK, SLM}. 

The method to identify coherent pairs used here is based on the Frobenius-Perron operator, the Ulam-Galerkin method and the thresholding method.
See \cite{FK}. 

In this work, we extend the concept of finite time coherent pairs to incorporate relative measure, and we call this relatively coherent pairs.  
This extension provides the theoretical framework to simply apply the definition of finite time coherence iteratively and at each stage we hierarchically define the relative measure on each element of the sub partitions which are developed. 
The results can be collected in a tree structure to emphasize the hierarchical nested nature of such partitions.

Also we present a adapted thresholding method that can identify the finite-time relatively coherent structures in successive scales of time-dependent systems. 
We study four examples in this paper, for both closed systems and open systems.
In the first and second example, 
we show our method identifies the fine scaled relatively coherent structures in a nonautonomous double gyres and a standard map. 
In the third example, an idealized zonal stratospheric flow, our method gives fine scaled details with respect to relative coherence. 
In our fourth example, we show the method is also efficient at an open system, corresponding to oceanographic flows in the Mexico Gulf.

This paper is organized as follows. In Sec. 2, we define relative coherent structures. 
Then we briefly describe the Frobenius-Perron operator, and the Ulam-Galerkin matrix estimate and the thresholding method. 
In Sec. 3, we give the details of the algorithm and a successive flow chart for explanation. 
From Sec. 4 to Sec. 7, we apply the method to four examples. 
Conclusions are given in Sec. 8. In Sec. 9, the appendix, we analyze the computational complexity.

%%%%%%%%%%%%%%%%%%%%%%%%%%%%%%%%%%%%%%%%%%%%%%%%%%%%%%%%%%%%%%%%%%%%%%%%%%%%%%%%%
%                                                                                                                      2 Relatively coherent structures
%%%%%%%%%%%%%%%%%%%%%%%%%%%%%%%%%%%%%%%%%%%%%%%%%%%%%%%%%%%%%%%%%%%%%%%%%%%%%%%%%

\section{Relatively coherent structures}

We define relatively coherent structures by specializing the definition of coherent pairs, which we now review \cite{FSM}.
Let ($\Omega$, $\mathcal{A}$, $\mu$) be a measure space, where $\mathcal{A}$ is a $\sigma$-algebra and $\mu$ is a normalized measure that is not necessarily invariant. 
The key to the specialization to relatively coherent pairs is the use of an iteratively defined relative measure, refined from the initial globe measure $\mu$.

Generally, we assume that $\Omega \subset \mathbb{R}^d$.
Given a time-dependent flow $\Phi (z,t; \tau ) : \Omega \times \mathbb{R} \times \mathbb{R} \to \Omega$, through the time epoch $\tau$ of an initial point $z$ at time $t$. 
A coherent pairs $(A_t, A_{t+\tau})$ can be considered as a pair of subsets of $\Omega$ such that,
$$\Phi (A_t,t; \tau ) \approx A_{t+\tau}.$$

\begin{definition}
$\cite{FSM} (A_t, A_{t+\tau})$ is a $(\rho _0, t, \tau)$-coherent pair if
{\begin{eqnarray}
\rho_\mu(A_t ,A_{t+\tau}):=\mu (A_t \cap \Phi(A_{t+\tau},t+\tau, -\tau))/\mu(A_t) \ge \rho_0
\label{CP}
\end{eqnarray}}
where the pair $(A_t, A_{t+\tau})$ are `robust ' to small perturbation and $\mu (A_t)= \mu (A_{t+ \tau}).$
\end{definition}

Note that the definition centrally depends on the full measure $\mu$ on $\omega$, and we will substitute successive relative measures on refinements.
Now we consider a relative measure on $\omega$ induced by $\mu$, where $\omega$ is a nonempty measurable subset of $ \Omega $. 
In this way enters refinements of the initial partition on successive scales, a {\bf relative measure} of $w$ to $\Omega $ is 
{\begin{eqnarray}
\mu _{\omega } (A):=\frac{\mu(A \cap \omega )}{\mu (\omega)}
\label{RM}
\end{eqnarray}}
 for all $A \in \mathcal{A}$.

From the above definition, it follows that the space $ (\omega, \mathcal{A} |_{\omega}, \mu_{\omega}) $  is also a measure space, 
where $\mathcal{A} |_{\omega}$ is the restriction of $\omega$ to $\mathcal{A}$ and $\mu_{\omega}$ is a normalized measure on $\omega$. 
We call the space $ (\omega, \mathcal{A} |_{\omega}, \mu_{\omega}) $, the $ relative \  measure \  space$.
Now, we define the $relatively \  coherent \   pairs$.

\begin{definition}
Relatively coherent structures are those $(\rho _0, t, \tau)$-coherent pairs defined Definition 2.1, 
with respect to given relative measures on a subset $\omega \subset \Omega$, of a given scale.
\end{definition}

To find relatively coherent structures in time-dependent dynamical systems, the basic tool is the Frobenius-Perron operator.
Let $(\Omega, \mathcal{A}, \mu)$ be a measure space and $\mu$ is a normalized Lebesgue measure. 
If $S: \Omega \to \Omega$ is a nonsingular transformation such that $\mu (S^{-1}(A))=0$ for all $A \in \mathcal{A}$ satisfying $\mu (A)=0,$ 
the unique operator $P: L^1 (\Omega) \to L^1 (\Omega)$ defined by,
{\begin{eqnarray}
\int _A Pf(x) \mu (dx) = \int _{S^{-1}(A)} f(x) \mu (dx)
\end{eqnarray}}
for all $A \in  \mathcal{A}$ is called the Frobenius-Perron operator corresponding to S, where $f(x) \in L^1 (\Omega)$ is a probability density function. 
See \cite{LM}.
In our case, $S$ can be considered as the flow map $\Phi$ and the formula above can be written as 
{\begin{eqnarray}
P_{t, \tau} f(z):= f(S^{-1}(z)) \cdot |det \  D(S^{-1}(z))|=f(\Phi(z,t+\tau; -\tau)) \cdot|det \  D\Phi(z,t+\tau; -\tau) |.
\end{eqnarray}}
Suppose X is a subset of M, let $Y$ be a set that includes $S(X)$. We develop partitions for $X$ and $Y$ respectively. 
In other words, let $\{ B_i \}_{i=1}^m$ be a partition for $X$ and $\{ C_j \}_{j=1}^n$ be a partition for $Y$.
The Ulam-Galerkin matrix follows a well-known finite-rank approximation of the Frobenius-Perron operator, which is of the form
{\begin{eqnarray}
P_{i,j}=\frac{\mu (B_i \cap S^{-1}(C_j))}{\mu (B_i)}
\end{eqnarray}}
where $\mu$ is the normalized Lebesgue measure on $\Omega$.
As usual, we numerically approximate $A_{i,j}$ by,
{\begin{eqnarray}
P_{i,j}=\frac{\#\{x_k: x_k \in B_i \   \& \   S(x_k) \in C_j\}}{\# \{x_k: x_k \in B_i\}}
\end{eqnarray}}
where the sequence $\{x_k\}$ is a set of test points (passive tracers). 
See \cite{DYZ}. \\

The following thresholding method from \cite{FSM} finds optimally coherent pairs in a time-dependent dynamical system, with respect to the chosen measure.
We will iteratively adapt the method to the relative measure.
This algorithm thresholds to the singular values and singular vectors of the matrix P obtained by the Ulam-Galerkin method : \\

\begin{enumerate}[i.]

\item Calculate the second singular value and corresponding left and right singular vectors $\{x_i\}$ and $\{y_j\}$ of the Ulam-Galerkin matrix.                 \\

\item Find values $\{(b_k, c_k)\}$ as pairs such that, 
{\begin{eqnarray}
\rho(X(b_k),Y(c_k))=\frac{\sum_{i:x_i>b_k\& j:y_j>c_k}p_i P_{ij}}{\sum_{i:x_i>b_k}p_i}
\label{THE1}
\end{eqnarray}}

by thresholding, where 
{\begin{eqnarray}
p_i=\mu(B_i), X(b_k)=\cup_{i:x_i>b_k} B_i \  \mbox{and}  \ Y(c_k)=\cup_{j:y_j>c_k} C_j.
\label{THE2}
\end{eqnarray}}

\item Choose a partition related to a pair $(b^*,c^*)$ of $(b_k,c_k)$ such that 
{\begin{eqnarray}
\rho^*=\max_{k}\{\ \rho(X(b_k),Y(c_k))\}.
\label{THE3}
\end{eqnarray}}
The partition is maximally coherent with respect to $\mu$ on $\Omega$ and the test set $\{(b_k, c_k)\}$.
\end{enumerate}

%%%%%%%%%%%%%%%%%%%%%%%%%%%%%%%%%%%%%%%%%%%%%%%%%%%%%%%%%%%%%%%%%%%%%%%%%%%%%%%%%
%                                                                                                                                3 Algorithm
%%%%%%%%%%%%%%%%%%%%%%%%%%%%%%%%%%%%%%%%%%%%%%%%%%%%%%%%%%%%%%%%%%%%%%%%%%%%%%%%%

\section{Algorithm}
We now describe how to find relatively coherent pairs. 
By the thresholding Eqs. \ref{THE1} -- \ref{THE3}, we have obtained optimal coherent pairs, which are defined as $(X_1,Y_1)$ and $(X_2,Y_2).$
$Y_1$ can be considered as the image of $X_1$ under a flow $\Phi$ in time-$\tau$. 
In order to find relatively coherent structures in $X_1,Y_1$, $X_2$ and $Y_2$, we define relative measures on each of these sets. 
Define relative measures $\mu_{X_1} (S)$ and $v_{Y_1}(T)$  on the coherent pair $(X_1, Y_1)$, according to Eq. \ref{RM}.
Then we have the measure spaces $ (X_1, \mathcal{A} |_{X_1}, \mu _{X_1}) $ and $ (Y_1, \mathcal{A} |_{Y_1}, v _{Y_1}) $ with $\mu_{X_1}$ and $v_{Y_1}$ the normalized measures
\footnote {The probability measure $v$ can be considered as the discretized image of $\mu$. The detail of construction of $v$ can be found in \cite{FSM}}.

The relative measures for both $X_1$ and $Y_1$ allow the adaptation of the thresholding methods on $ (X_1, \mathcal{A} |_{X_1}, \mu _{X_1})$  
and $ (Y_1, \mathcal{A} |_{Y_1}, v _{Y_1}) $ 
under the relative measures. 
Then follows two relatively coherent pairs in the previous coherent pair $(X_1, Y_1)$, which are now named as $(X_{11}, Y_{11})$ and $(X_{12}, Y_{12})$. 
Also, the coherent pair $(X_2, Y_2)$ can be divided to two relatively coherent pairs, $(X_{21}, Y_{21})$ and $(X_{22}, Y_{22})$.

Next, we repeat the building process above, 
that is, to define eight normalized measures on each of  $X_{11}$, $Y_{11}$, $X_{12}$, $Y_{12}$, $X_{21}$, $Y_{21}$, $X_{22}$ and $Y_{22}$, respectively, 
such that they become new spaces with corresponding relative measures, 
we then apply the adapted thresholding method on these new spaces to get more relatively coherent structures.

Now, we state our hierarchical adaptation method as an algorithm, and for convenience, we use $(X_i,Y_i)$ and $(X_j,Y_j)$ to denote two coherent pairs,
i and j which can be stated to emphasize the hierarchy tree. \\

{\bf Algorithm 1}

\begin{figure}[htb]
  \centering
  \subfloat{\label{fig:gull}\includegraphics[width=1.0\textwidth]{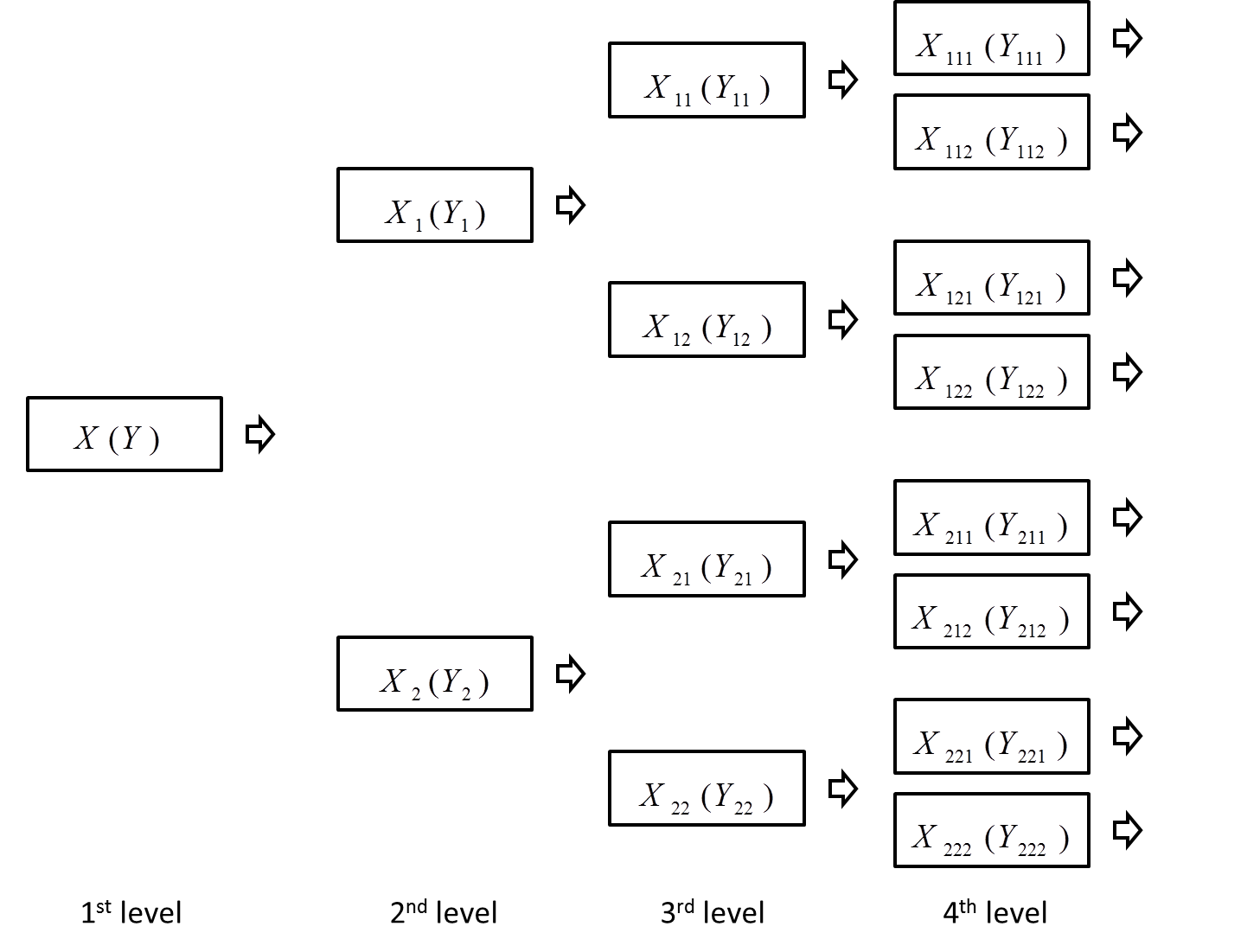}}         
\caption{Successive relatively coherent sets tree as per Algorithm 1.}
  \label{FlowChart}
\end{figure}

\begin{enumerate}
\item Define relative measures $\mu_{X_i}$ and $v_{Y_i}$ 
and relative measure spaces $ (X_i, \mathcal{A} |_{X_i}, \mu _{X_i}) $ and $ (Y_i, \mathcal{A} |_{Y_i}, v _{Y_i})$ for $(X_i,Y_i)$, where

{\begin{eqnarray}
\mu_{X_i} (S)=\frac{\mu(S)}{\mu (X_i)} \   \mbox{and} \   v_{Y_i}(T)=\frac {v(T)}{v(Y_i)}
\end{eqnarray}}

for all $S \subset X_i$ and for all $T \subset Y_i$.

\item Apply the adapted thresholding method on $ (X_i, \mathcal{A} |_{X_i}, \mu _{X_i}) $ and $ (Y_i, \mathcal{A} |_{Y_i}, v _{Y_i})$ for $(X_i,Y_i)$ 
to develop refined coherent pairs $(X_{ii},Y_{ii})$ and $(X_{ij},Y_{ij})$.

\item Repeat the above two steps for the pair $(X_j,Y_j)$ to obtain further refined coherent  pairs $(X_{ji},Y_{ji})$ and $(X_{jj},Y_{jj})$.
 
According to the first three steps, a relatively coherent structure can be denoted by
{\begin{eqnarray}
(X_{k_1 k_2 ... k_q}, Y_{k_1 k_2 ... k_q})
\label{KQ}
\end{eqnarray}}
and $k_p \in \{ i, j \}, p=1, 2, ...q$, through $q$-steps of the algorithm, $q+1$ levels deep into the tree. 
See Figure \ref{FlowChart}.
That is, the subscript $k_1 k_2 ... k_q$ can be any possible finite $q$ permutations of $i$ and $j$. 
We usually choose $i=1$ and $j=2$. 
Figure \ref{FlowChart} is a flow chart depicting four levels, which describes the steps to find relatively coherent pairs in finer `scales'. 
The chart emphasizes a `tree' structure.

However, we cannot repeat the algorithm forever, so we must decide a stopping criterion. 
%And also, sometimes, we need to ignore those relatively coherent pairs we are not interested in in order to make the algorithm more efficient.
The following step is as a completion for the algorithm.

\item Stop a given branch if in Eq. \ref{THE3},

{\begin{eqnarray}
\rho^*=\max_{k}\{\ \rho(X_{k_1 k_2 ... k_q}(b_k), Y_{k_1 k_2 ... k_q}(c_k))\}.
\end{eqnarray}}
 
is such that,

{\begin{eqnarray}
\rho^*<\rho _0.
\label{Thresh}
\end{eqnarray}}

where $\rho_0 \in (0, 1)$ is a threshold from Eq. \ref{CP}, descriptive of optimal coherence which is not very coherent.

%{\begin{eqnarray}
%\rho^*<\rho _0 \  \mbox{or} \  \rho^*>\tilde{\rho}.
%\label{Thresh}
%\end{eqnarray}}
%where $\rho_0 \in (0, 1)$ is from Eq. \ref{CP} and $\tilde{\rho} \in (\rho _0 ,1]$. 
%$\tilde{\rho}$ is chosen as a threshold to exclude those relatively coherent pairs with high coherent values from Eq. \ref{THE3}.
%Then we can focus on the regions we specified by Eq. \ref{Thresh}.

\end{enumerate}
The Gulf Example in Figure 7 on how such stopping criterion leads to the number of coherent pairs less than $2^q.$

%%%%%%%%%%%%%%%%%%%%%%%%%%%%%%%%%%%%%%%%%%%%%%%%%%%%%%%%%%%%%%%%%%%%%%%%%%%%%%%%%
%                                                                                                                     Example 1
%%%%%%%%%%%%%%%%%%%%%%%%%%%%%%%%%%%%%%%%%%%%%%%%%%%%%%%%%%%%%%%%%%%%%%%%%%%%%%%%%

\section{Example 1. The Nonautonomous Double Gyre.}

Consider the nonautonomous double gyre system,
\begin{eqnarray}
&& \dot{x}=- \pi A \sin ( \pi f(x,t)) \cos(\pi y) \nonumber \\
&&\dot{y}=\pi A \cos ( \pi f(x,t)) \sin(\pi y) 
\end{eqnarray}
where $f(x,t)= \epsilon \sin(\omega t) x^2 + (1-2 \epsilon \sin (\omega t)) x$, $\epsilon = 0.25$, $\omega=2 \pi$ and $A=0.25$.
See \cite{FK}. 

\begin{figure}[htb!]
  \centering

  \subfloat[t=0]{\label{fig:gull}\includegraphics[width=1\textwidth]{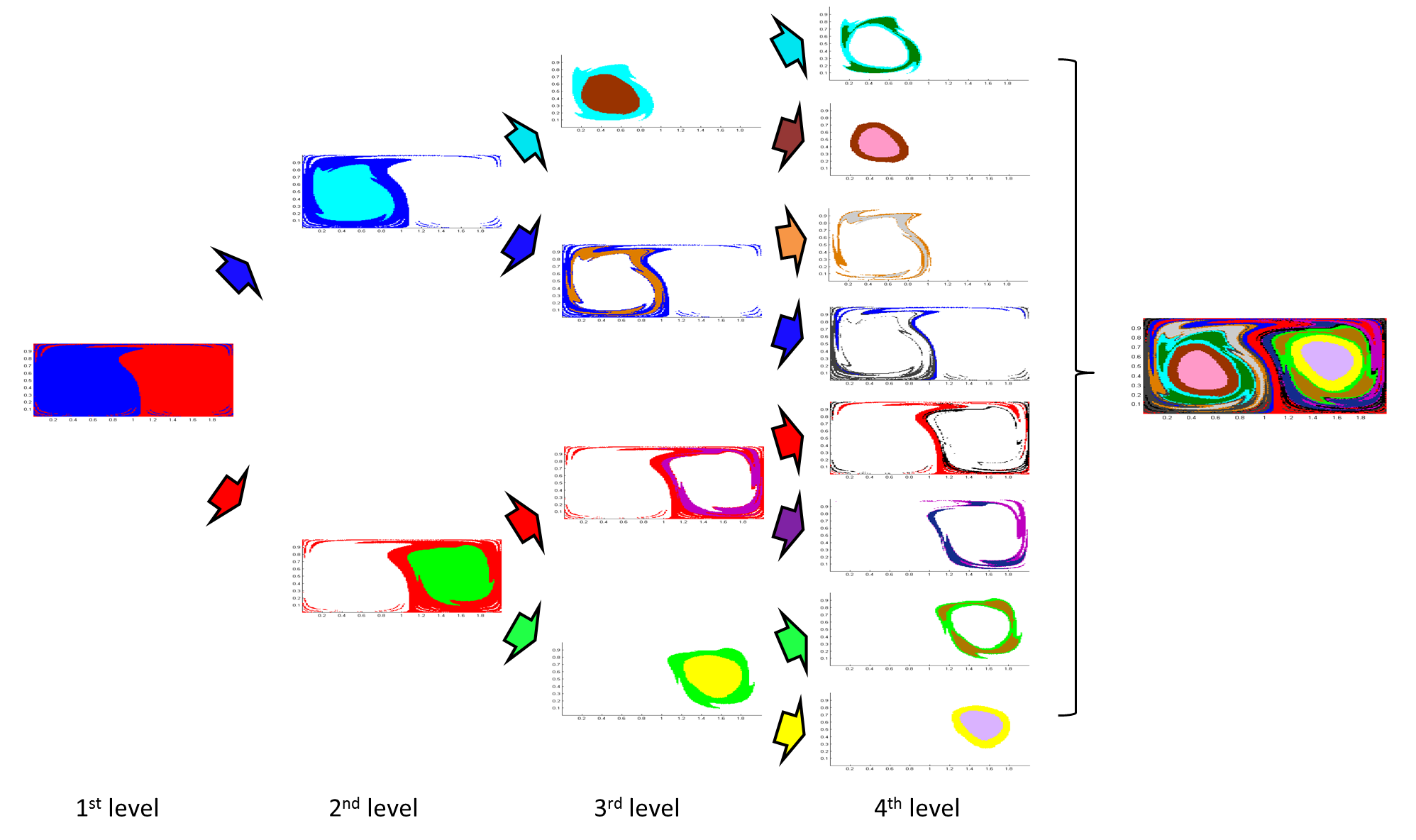}}    \\            
  \subfloat[t=10]{\label{fig:tiger}\includegraphics[width=1\textwidth]{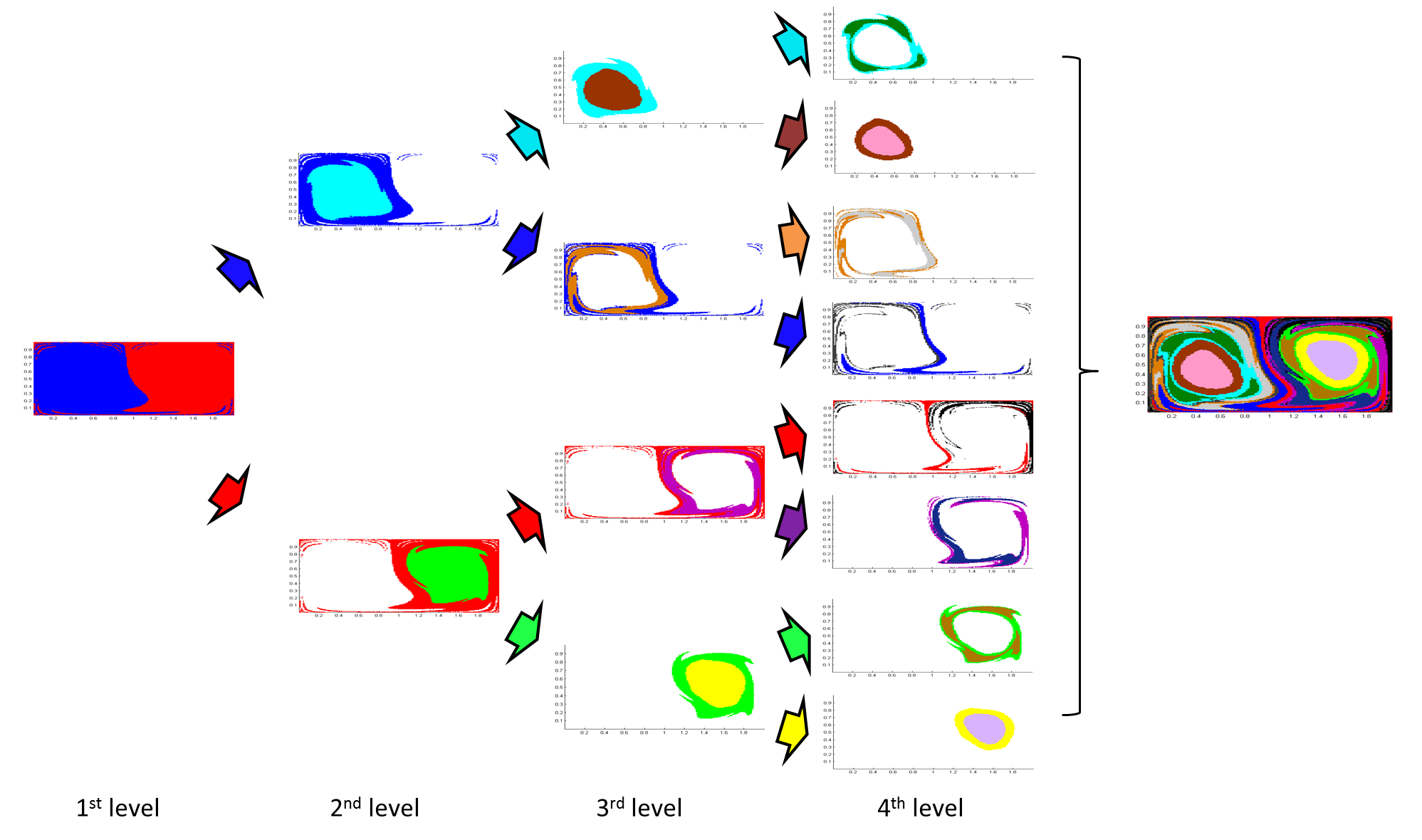}} 

  \caption{The upper figure is the initial status of the double gyre and the lower figure is the image under the time-$\tau$ flow, where $\tau=10$.  The same color areas between the two associates one relatively coherent pair. In this case, we have 16 different relative coherent structures with different colors. By following the colored arrows, we can see the relative coherent structures through four levels of refinement.}
  \label{DG}
\end{figure}

Let the initial time be $t=0$ and the final time be $t=10$. 
We use 24,200 identical triangles $\{B_i \}_{i=1}^{24200}$ to cover the region $[0,2] \times [0,1]$ by Delaunay triangulation. 
Then randomly and uniformly we choose 10,000,000 points in the region as our initial conditions.
For good sampling, the relationships between the given grid and the necessary number of points we choose is discussed in Appendix I, on computational complexity.
We numerically calculate the final status of these points by the Runge-Kutta method to estimate the flow. 
Because the double gyre model is an area-preserving system, the same triangulation can be used as the image partition $\{C_j \}_{j=1}^{24200}$.
The Ulam-Galerkin's transition matrix estimates the Frobenius-Perron operator that has the size $24200$ by $24200$. 
We apply the thresholding method on the matrix to find the first two coherent pairs in the initial status and final status. 
See Figure \ref{DG}. 
The first two coherent pairs are colored blue and red in the first level of both the upper and lower charts. 
We define the left part of the initial status as $X_1$ and the left part of the final status as $Y_1$, both of which are filled with blue as halves of the initial status and final status. 
Thus, two relative measures can be defined on each of these two parts separately, 

\begin{eqnarray}
 \mu_{X_1} (S)=\frac{\mu(S)}{\mu (X_1)}, \   \   \  v_{Y_1}(T)=\frac {v(T)}{v(Y_1)}.  
\end{eqnarray}

On the other hand, we can do the same in the right red-filled regions which we call the initial status $X_2$ and  the final status $Y_2$ to develop another two relative measures as follow,
 
\begin{eqnarray}
 \mu_{X_2} (S)=\frac{\mu(S)}{\mu (X_2)},  \   \   \    v_{Y_2}(T)=\frac {v(T)}{v(Y_2)}.
\end{eqnarray}

Now $X_1$ and $Y_1$ can be considered as the initial status and the final status of a refined relative sub-`system'. 
Likewise for $X_2$ and $Y_2$. 
By the same process as with the whole double gyres system, we can get some new coherent structures in the new system consisting of $X_1$ and $Y_1$. 
In Figure \ref{DG}, following the first blue arrow between first level and second level of both flow charts, $X_1$, the blue half on the first level of the upper chart is divided by blue and light blue, we define the blue part as $X_{11}$ and the light blue part as $X_{12}$ in the second level. 
Correspondingly, we have $Y_{11}$ which is blue and $Y_{12}$ which is light blue in the second level of the lower chart of Figure \ref{DG}  as the outcome states of $X_{11}$ and $X_{12}$. 
$X_{11}$ and $Y_{11}$ are relatively coherent structures, 
and so are $X_{12}$ and $Y_{12}$.

As above, we can develop an $X_{21}$ that is red and an $X_{22}$ that is light green  from ${X_2}$ in the second level of the upper chart of Figure \ref{DG}; 
and  $Y_{21}$ that is red and $Y_{22}$ that is light green  from ${Y_2}$ in the second level of the lower chart of Figure \ref{DG}. The same subscript means $X_{21}$ and $Y_{21}$ are a relatively coherent pair, so are $X_{22}$ and $Y_{22}$. Now we have four relatively coherent pairs in the second level.

We can eventually get the tree structures in Figure \ref{DG} by repetition of the process. 
There are 8 relatively coherent structures shown in the third level with different colors and 16 relatively coherent structures in the fourth level with different colors. 
In Figure \ref{DG}, we can see the egg shaped relatively coherent structures which are four resonance `islands' as expected in such hamiltanian twist maps, \cite {M}.
Even finer structures will be revealed by further refinement and sufficient sampling to allow appropriate resolution. 
Appropriate sampling in a given refinement scale is discussed in terms of computational complexity in developing a given Ulam-Galerkin's matrix for a given fine grid.

%%%%%%%%%%%%%%%%%%%%%%%%%%%%%%%%%%%%%%%%%%%%%%%%%%%%%%%%%%%%%%%%%%%%%%%%%%%%%%%%%
%                                                                                                                     Example 2
%%%%%%%%%%%%%%%%%%%%%%%%%%%%%%%%%%%%%%%%%%%%%%%%%%%%%%%%%%%%%%%%%%%%%%%%%%%%%%%%%

\section{Example 2. The Standard Map.}

Consider the standard map, 
\begin{eqnarray}
&& p_{n+1}=p_n +K \sin (\theta _n) \nonumber \\
&&\theta _{n+1}= \theta _n + p_{n+1}
\label{ST}
\end{eqnarray}
where $p_n$ and $\theta _n$ are taken modulo $2 \pi$. 
See \cite{M}, and this is a map on the torus $\Omega=[0,1) \times [0,1).$
We will study the case that $K=1.2$ in our example, 
as this is well known to be $K=1.2 >K_{cr}=0.971635$ shortly after the breakup of the last `golden' torus allowing momentum boosting orbits and a mixed chaotic and ordered phase space including periodic elliptic islands, \cite{M}.

We choose 10,000,000 points, randomly and uniformly in the area $[0,1) \times [0,1)$. 
We set $X:=[0,1) \times [0,1)$ and $Y:=[0,1) \times [0,1)$ and then uniformly triangulate to cover both the initial status and the mapping image. 
Therefore, enumerating the triangulation grid $\{B_i \}_{i=1}^{20000}$ and $\{C_j \}_{j=1}^{20000}$,
let the initial points do 10 iterations, after which we obtain the final status of the system and the Ulam-Galerkin's matrix is 20,000 by 20,000.

\begin{figure}[htb]
 \centering

 \subfloat[At the beginning]{\label{fig:gull}\includegraphics[width=.7 \textwidth]{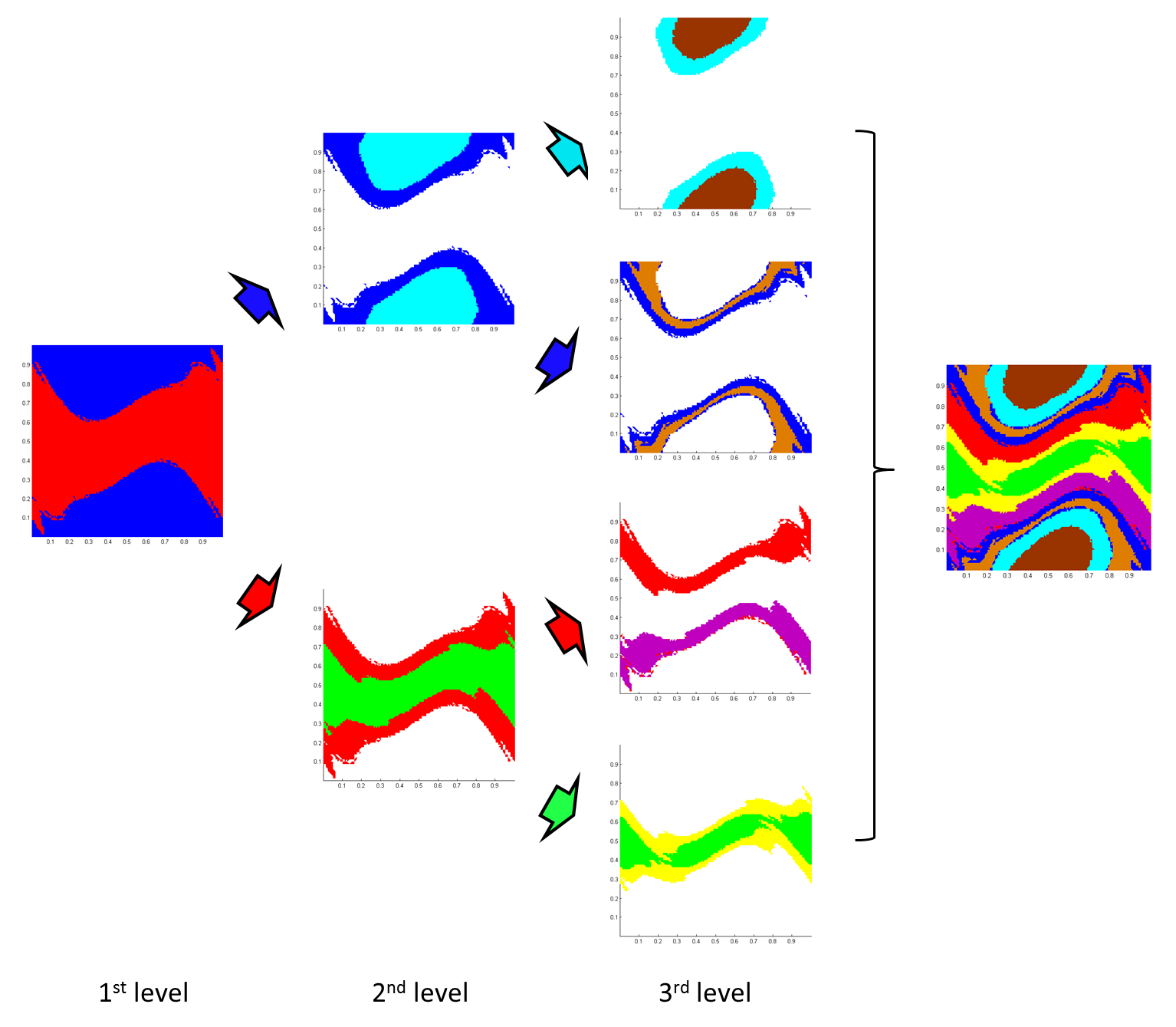}}       \\
 \subfloat[After 10 iterations]{\label{fig:mouse}\includegraphics[width=.7 \textwidth]{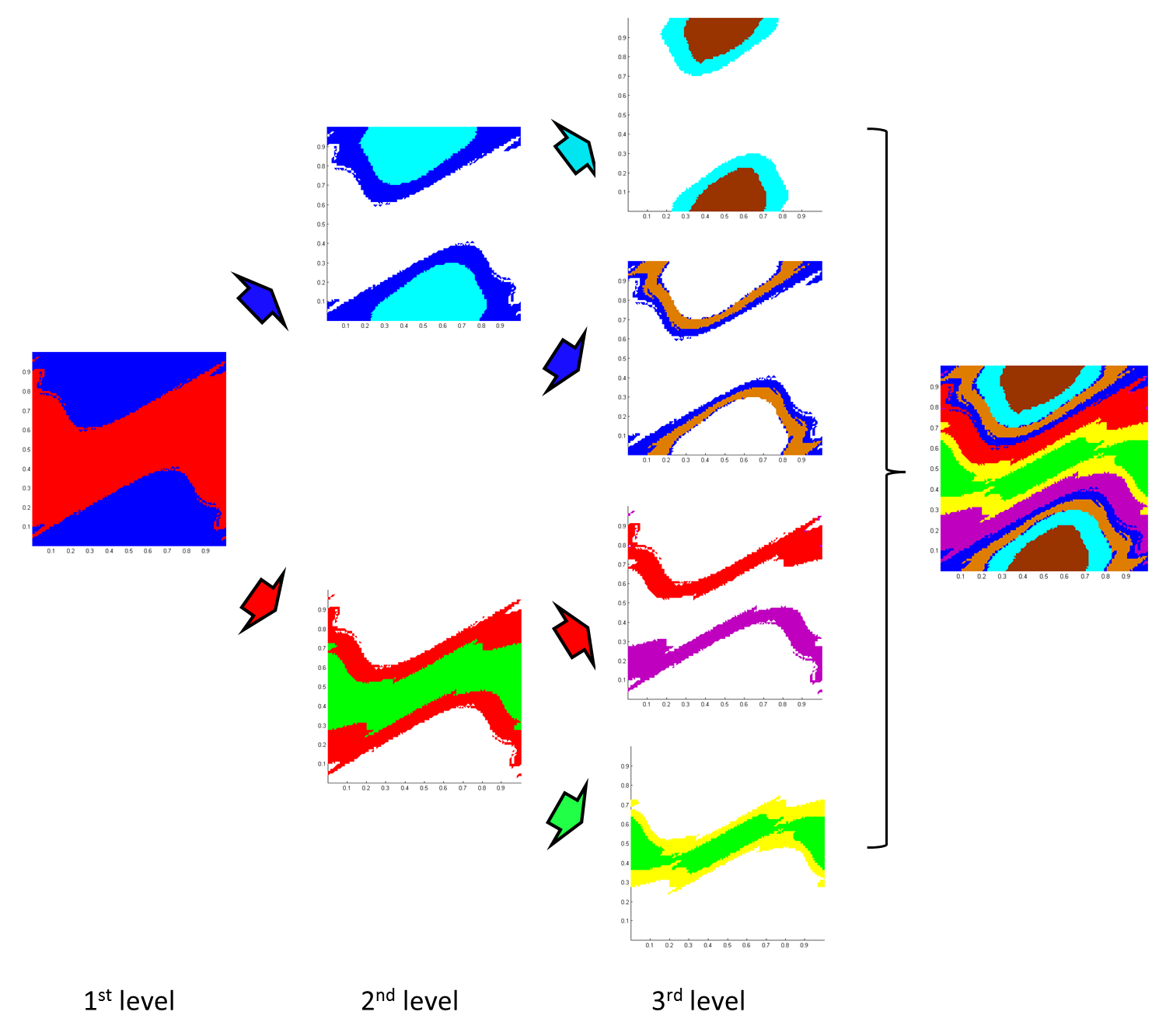}}

  \caption{Standard Map Eq. \ref{ST}, structured relative coherence hierarchy tree arranged as in (a) - (b). }
  \label{SM1}
\end{figure}

The first level of Figure \ref{SM1} shows two large coherent pairs in the Standard Map, which are filled by blue and red, respectively.
By prior knowledge on this benchmark problem, it is clear that the boundary between the primary blue and red coherence estimates the cantorus remnant of the golden area resonance, an expected act came as this is known to remain a primary pseudo-barrier to transport when $K=1.2$ still not much larger than $K_{cr}=0.971635.$

Iteratively repeating the process according to the algorithm for each of $X_1$, $Y_1$, $X_2$ and $Y_2$ yields 4 relatively coherent structures 
which are colored blue, light blue, red and green in the second level. 
The third level in Figure \ref{SM1} tells us there are 8 different such structures in total in the first three levels.
The outcome partition here shows a familiar depiction of the resonance layers known to be due to cantorus pseudo-barriers 
which cause the famously slow transport for the standard map, \cite{M}.

%%%%%%%%%%%%%%%%%%%%%%%%%%%%%%%%%%%%%%%%%%%%%%%%%%%%%%%%%%%%%%%%%%%%%%%%%%%%%%%%%
%                                                                                                                     Example 3
%%%%%%%%%%%%%%%%%%%%%%%%%%%%%%%%%%%%%%%%%%%%%%%%%%%%%%%%%%%%%%%%%%%%%%%%%%%%%%%%%

\section{Example 3. An Idealized Stratospheric Flow.}
 
Next, consider the Hamiltonian system
\begin{eqnarray}
dx/dt=-\partial \Phi / \partial y \nonumber \\
dy/dt=\partial \Phi / \partial x,
\end{eqnarray}
where
\begin{eqnarray}
\Phi(x,y,t)=c_3y-U_0Ltanh(y/L)+A_3U_0Lsech^2(y/L)cos(k_1x) \nonumber \\
 +A_2U_0Lsech^2(y/L)cos(k_2x-\sigma _2t) \nonumber \\
+A_1U_0Lsech^2(y/L)cos(k_1x-\sigma _1t)
\label{ROWA}
\end{eqnarray}

This is a quasiperiodic system that represents an idealized zonal stratospheric flow, \cite{FSM}. 
There are two known Rossby wave regimes in this system. 
Let $U_0=63.66, c_2=0.205U_0, c3=0.7U_0, A3=0.2, A2=0.4, A1=0.075$ and the other parameters in Eq. \ref{ROWA}  be the same as stated in \cite{RBB}. 

\begin{figure}[htb!]
  \centering

  \subfloat[t=0]{\label{fig:gull}\includegraphics[width=0.6\textwidth]{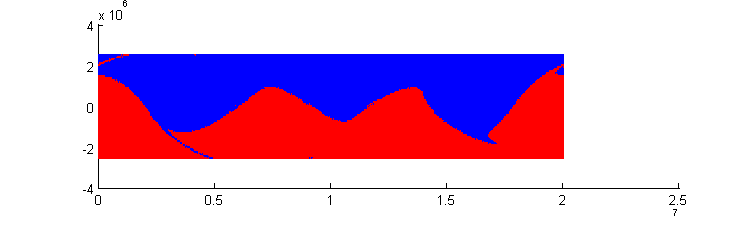}}                
  \subfloat[t=0]{\label{fig:tiger}\includegraphics[width=0.6 \textwidth]{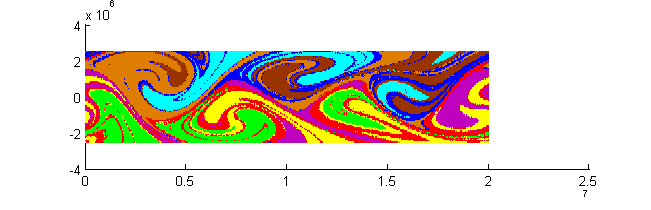}} \\

  \subfloat[t=10 \ days]{\label{fig:mouse}\includegraphics[width=0.6\textwidth]{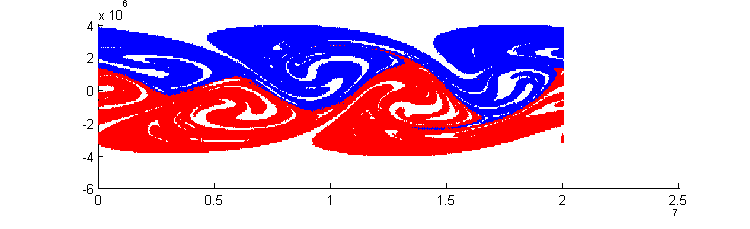}}
   \subfloat[t=10 \ days]{\label{fig:mouse}\includegraphics[width=0.6\textwidth]{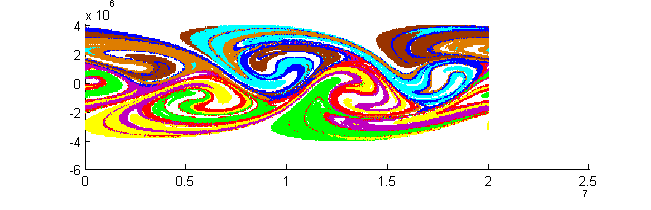}}
  \caption{Relative coherence in the Rossby system Eq. \ref{ROWA}. The fist level partition (left hand side) and third level partition (right hand side) of both the initial and final status of the zonal flow. Compare to the hierarchical structure as emphasized in Figure \ref{RW2}.}
  \label{RW1}
\end{figure}

\begin{figure}[htb!]
  \centering

 \subfloat[t=0]{\label{fig:gull}\includegraphics[width=1.0\textwidth]{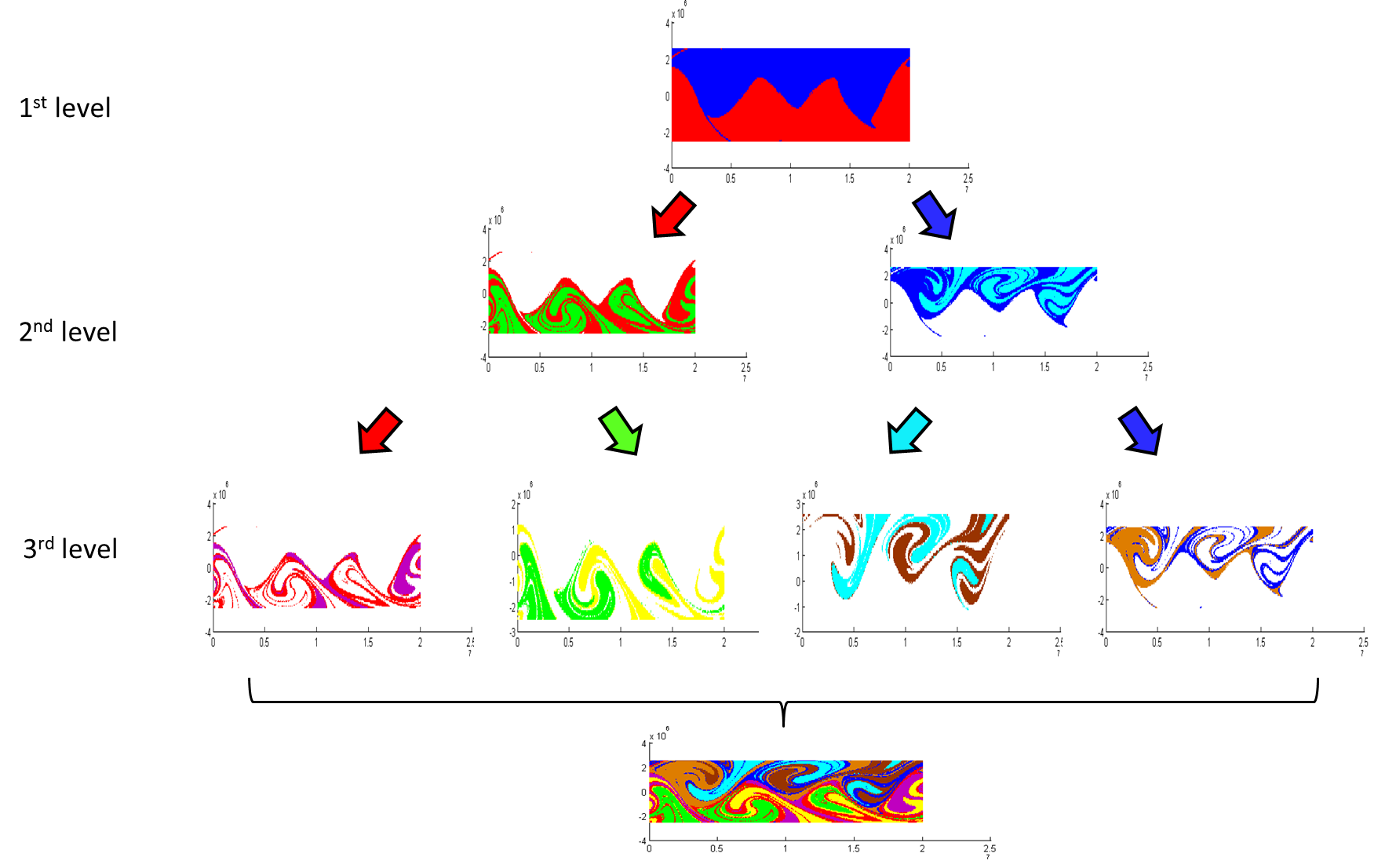}}       \\

 \subfloat[t=10 \ days]{\label{fig:mouse}\includegraphics[width=1.0\textwidth]{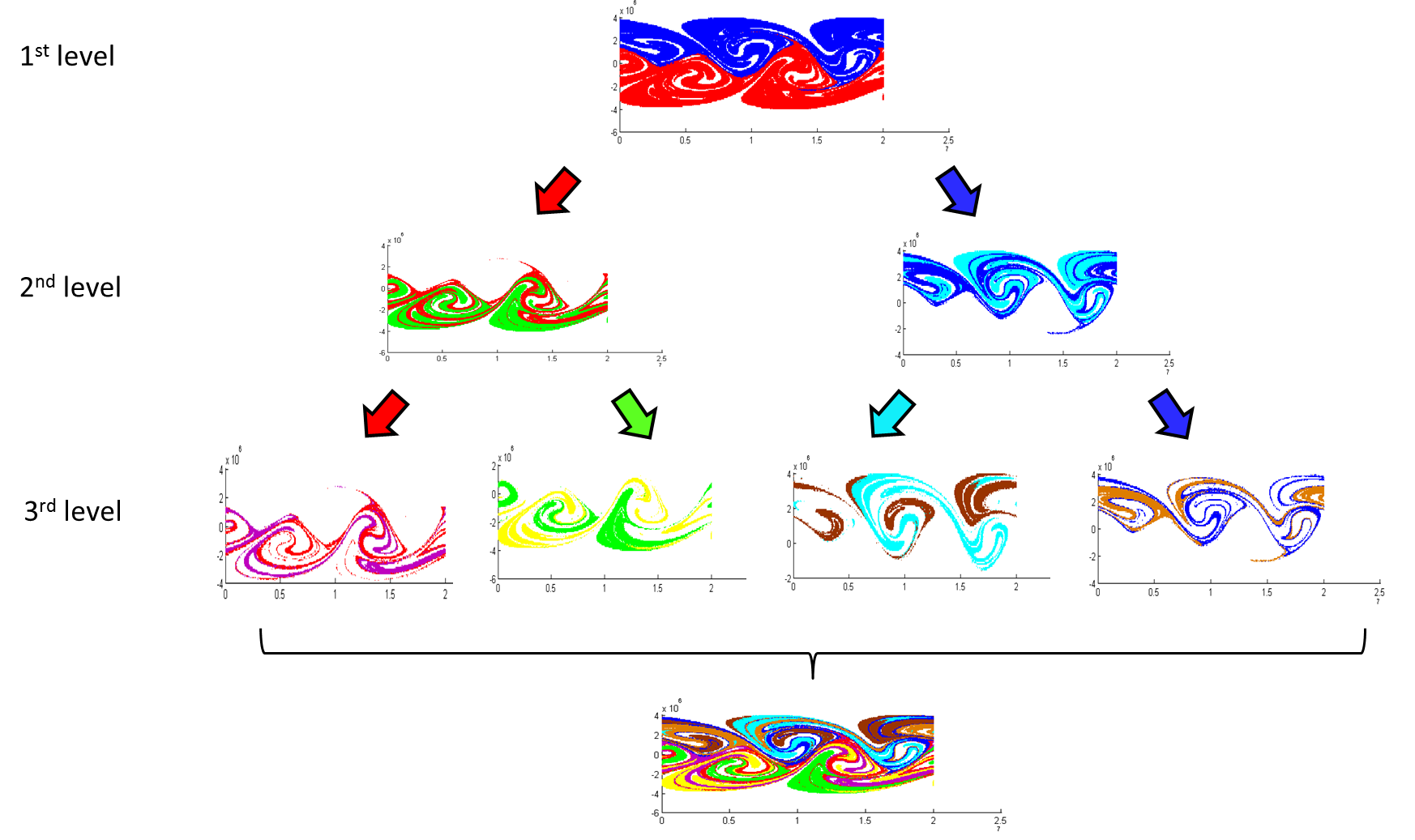}}

  \caption{Relative coherence hierarchy in the Rossby system Eq. \ref{ROWA}. Coloring and hierarchy tree structure as in (a)-(b). Compare also to Figure \ref{RW1}. }
  \label{RW2}
\end{figure}

We choose 20,000,000 points in the domain $X=[0, \  6.371 \pi*10^6] \times [-2.5*10^6, \  2.5*10^6 ]$ of the flow and use 32,640 triangles as the partition $\{B_i \}_{i=1}^{32640}$ for the initial status points and 39,694 triangles as the partition $\{C_j \}_{j=1}^{39694}$ for the final status of the points.
Note that this system is `open' relative to the windows $X$ chosen. 
The two coherent pairs are colored blue and red which are defined as $(X_1, Y_1)$ and $(X_2, Y_2)$  in the first level of Figure \ref{RW1}. 
Again, we now build the relative measures and tree of relatively coherent pairs. 
By applying the method as we have done with the previous two examples, we develop four and eight different coherent structures for the second level and the third level, respectively.
In Figure \ref{RW2}, we repeat the color scheme in both the initial status and final status in each of the second and third levels for the relative coherent structure.
Thus, we now see a much finer scaled relative coherence in the dynamical system than previously seen.

%%%%%%%%%%%%%%%%%%%%%%%%%%%%%%%%%%%%%%%%%%%%%%%%%%%%%%%%%%%%%%%%%%%%%%%%%%%%%%%%%
%                                                                                                                     Example 4
%%%%%%%%%%%%%%%%%%%%%%%%%%%%%%%%%%%%%%%%%%%%%%%%%%%%%%%%%%%%%%%%%%%%%%%%%%%%%%%%%

\section{Example 4. The Gulf of Mexico.}
\begin{figure}[htb!]
  \centering

  \subfloat[]{\label{fig:gull}\includegraphics[width=1.0\textwidth]{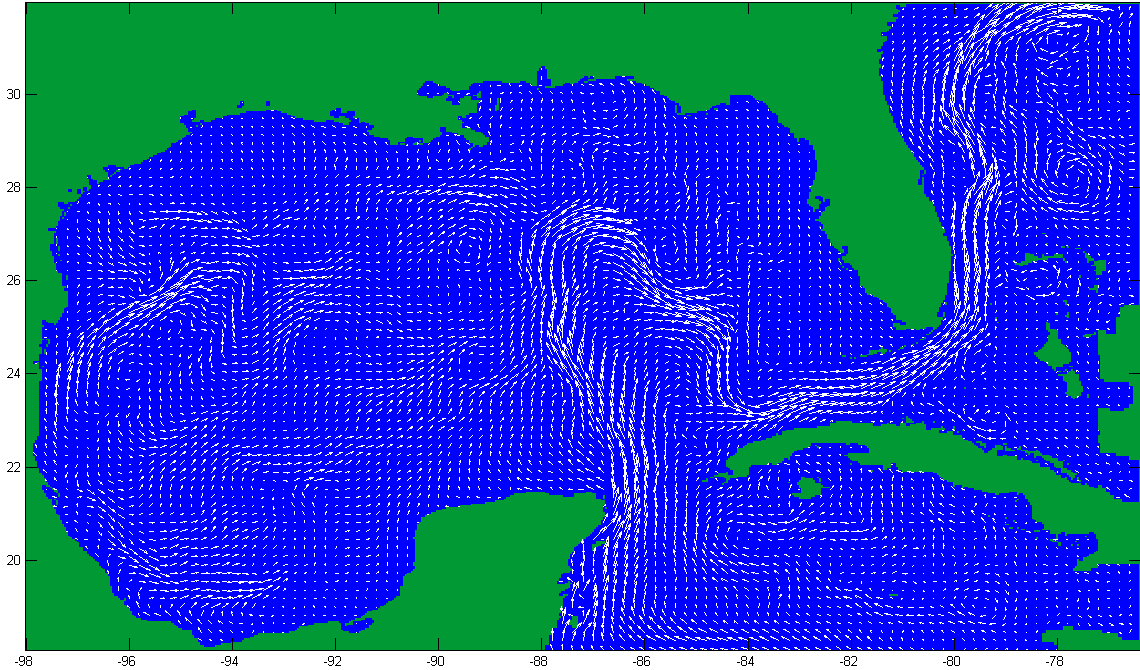}}    \\            

  \caption{Vector field describing surface flow in the Gulf of Mexico on May 24, 2010, computed using the HYCOM model [HYCOM, 2010]. Note the coherence of the Gulf Stream at this time. Oil spilling from south of Louisiana could flow directly into the Gulf Stream and out towards the Atlantic. This is an open system relative to this window shown. Horizontal and vertical units are degrees longitude (negative indicates west longitude) and degrees latitude (positive indicates north latitude), respectively.}
  \label{GulfVF}
\end{figure}

\begin{figure}[htb!]
  \centering

  \subfloat[t=0]{\label{fig:gull}\includegraphics[width=1.0\textwidth]{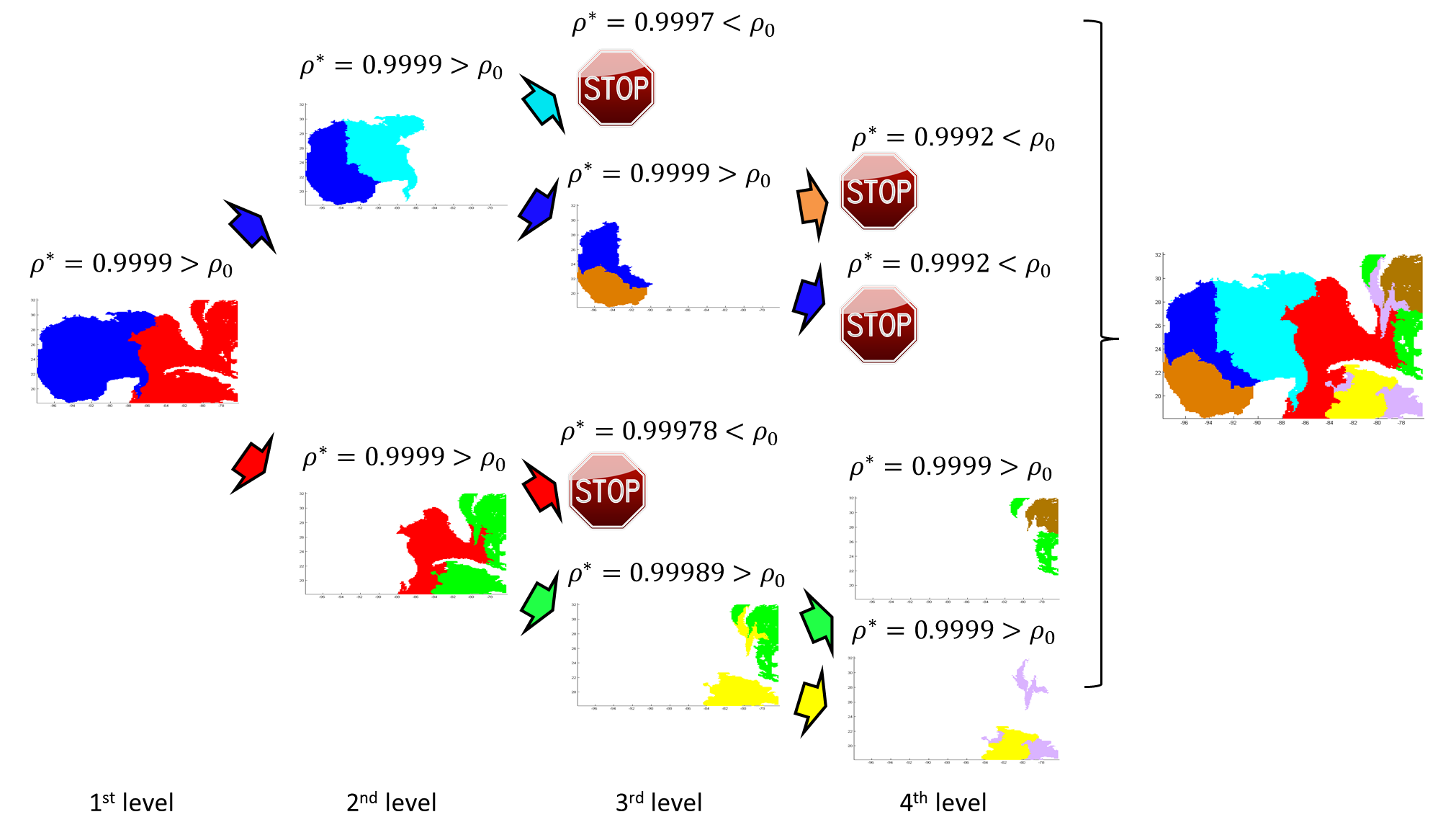}}    \\            
  \subfloat[t=6 \ days]{\label{fig:tiger}\includegraphics[width=1.0\textwidth]{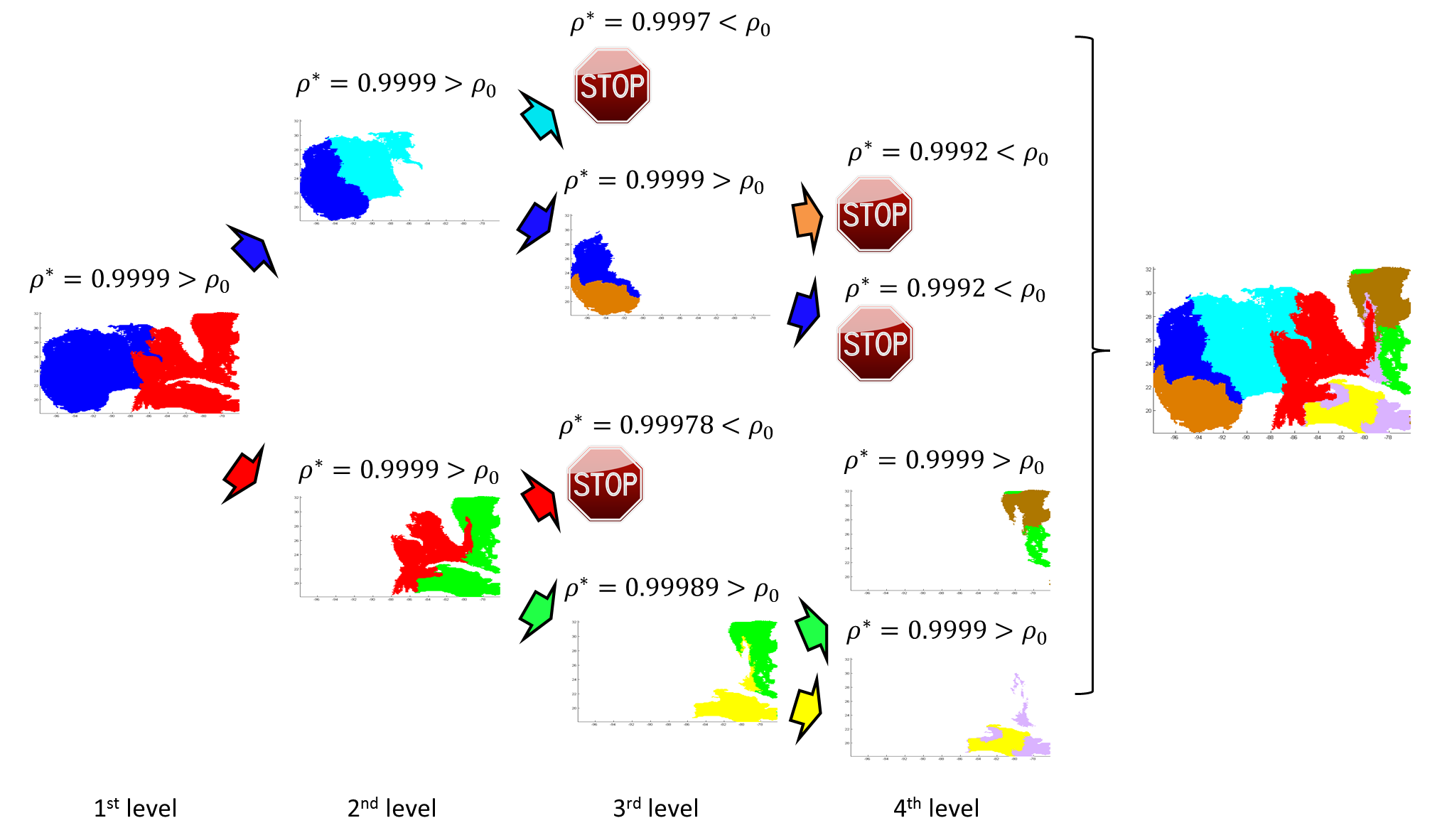}} 

  \caption{Hierarchical relative coherence in the Gulf of Mexico following the flow according to vector field data as illustrated in Figure \ref{GulfVF}. Tree structure and relatively coherent pairs coloring as in (a)-(b).}
  \label{Gulf}
\end{figure}

In our last example, we consider the Mexico Gulf. 
The data is the same as was used in \cite{BLSB} and formed by the method in \cite{B} and \cite{HA} .
The difference between the Gulf model and the above three examples is the Gulf is a open system, 
that is, there is water entering and exiting the region. 
See Figure \ref{GulfVF}.
This is the reason why in Figure \ref{Gulf}, the shapes of the whole Gulf water of the initial status and final status are slightly different at the bottom region and top region. 

We choose 20,000,000 points uniformly and randomly in the water region as the initial status, on more details of data. See \cite{BLSB}.
The final status is the positions of these points after 6 days.
We use 21,645 triangles $\{B_i \}_{i=1}^{32867}$ as a partition of $X$  and 21,210 triangles $\{C_j \}_{j=1}^{32359}$ as a partition of $Y$.
After applying our subdivision method on these triangles, the results are shown in Figure \ref{Gulf}. 
In this example, we set the $\rho_0=0.9998$ as the threshold the stopping criterion. 
Therefore, the number of relatively coherent pairs is not equal $2^q$, where $q$ is defined in Eq. \ref{KQ}.

%%%%%%%%%%%%%%%%%%%%%%%%%%%%%%%%%%%%%%%%%%%%%%%%%%%%%%%%%%%%%%%%%%%%%%%%%%%%%%%%%
%                                                                                                                     Conclusions
%%%%%%%%%%%%%%%%%%%%%%%%%%%%%%%%%%%%%%%%%%%%%%%%%%%%%%%%%%%%%%%%%%%%%%%%%%%%%%%%%

\section{Conclusions}
We have defined a concept of relative coherence based on relative measure, as a generalization of coherent pairs.
We also have introduced a recursive method of detecting relatively coherent structures under flows in a finite time, 
which is based on relative measures, with respect to the restricted Frobenius-Perron operator. 
Relative measures are used to build a hierarchy of relatively coherent pairs in successive levels, which can be illustrated in a natural tree structure of relative coherence.

We have demonstrated the method with the double gyres, the standard map, a Rossby wave system and data from the Gulf of Mexico. These examples have included hierarchical structure, open and closed systems, a system known only through data, and use of the stopping criterion.

%%%%%%%%%%%%%%%%%%%%%%%%%%%%%%%%%%%%%%%%%%%%%%%%%%%%%%%%%%%%%%%%%%%%%%%%%%%%%%%%%
%                                                                                                                     Appendix
%%%%%%%%%%%%%%%%%%%%%%%%%%%%%%%%%%%%%%%%%%%%%%%%%%%%%%%%%%%%%%%%%%%%%%%%%%%%%%%%%

\section{Appendix I: On Computational Complexity.}
In this section, we analyze the computational complexity to properly develop an Ulam-Galerkin matrix. 
The question is how many initial points should we use for a given grid. 
In other words, how many initial points is enough to well represent the whole domain for a given grid. 
Intuitively, we may wish to add `as many as possible'. 
Generally, finer grids require exponentially more points, depending on the dimensionality, 
but also on the local stretching of the map.
The following discussion is base on the Lipschitz constant and Gronwall's inequality. 
See \cite{P}.
There are some other works related to this topic. 
See \cite{DJ2, GDK, GK, J}.

Consider a grid that consists of square boxes with length $q$. 
The triangulation we use is handled similarly, but rectangles will simplify this discussion even if triangles allow for the powerful Delaunay triangulation algorithms in practice.
For convenience, we use boxes instead of triangles, but one box can be easily changed to two triangles by cutting through the diagonal. 

There are several steps to find how many points we need. 
Let $f: X \times [t_0,t] \to X$ be a flow and X be compact.

\begin{enumerate}[I.]

\item To find an $\varepsilon$ such that for $\forall x_1, x_2 \in X, |x_1(t_0)-x_2(t_0)|< \varepsilon$ we have 
{\begin{eqnarray}
|x_1(t)-x_2(t)| < q,
\end{eqnarray}}
where $x_1(t_0)$ and $x_2(t_0)$ are the initial status of $x_1, x_2$ and $x_1(t), x_2(t)$ are the final status of $x_1, x_2$.
The purpose of this step is to avoid that the distance $d$ between the final positions of two closed initial points under the flow change dramatically, in our case, it means $d \geq q$. 
See Figure \ref{PA1} case (a). 
So, we must control $\varepsilon$ such that case (b) holds for all initial points.

\begin{figure}[htb!]
  \centering

  \subfloat[$d>q$. A failure scenario]{\label{fig:gull}\includegraphics[width=.5\textwidth]{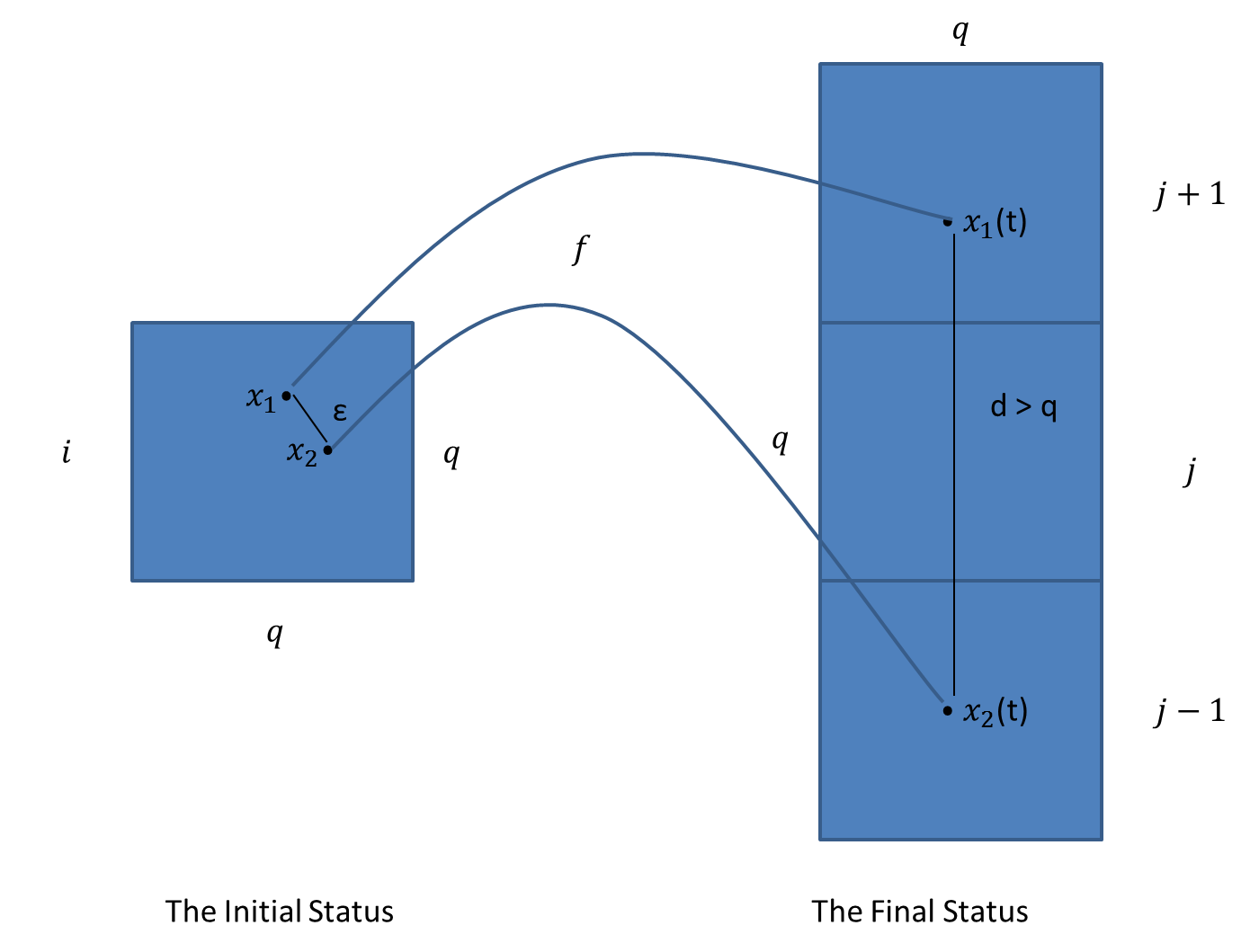}}             
  \subfloat[ $d<q$. The success scenario]{\label{fig:gull}\includegraphics[width=.5\textwidth]{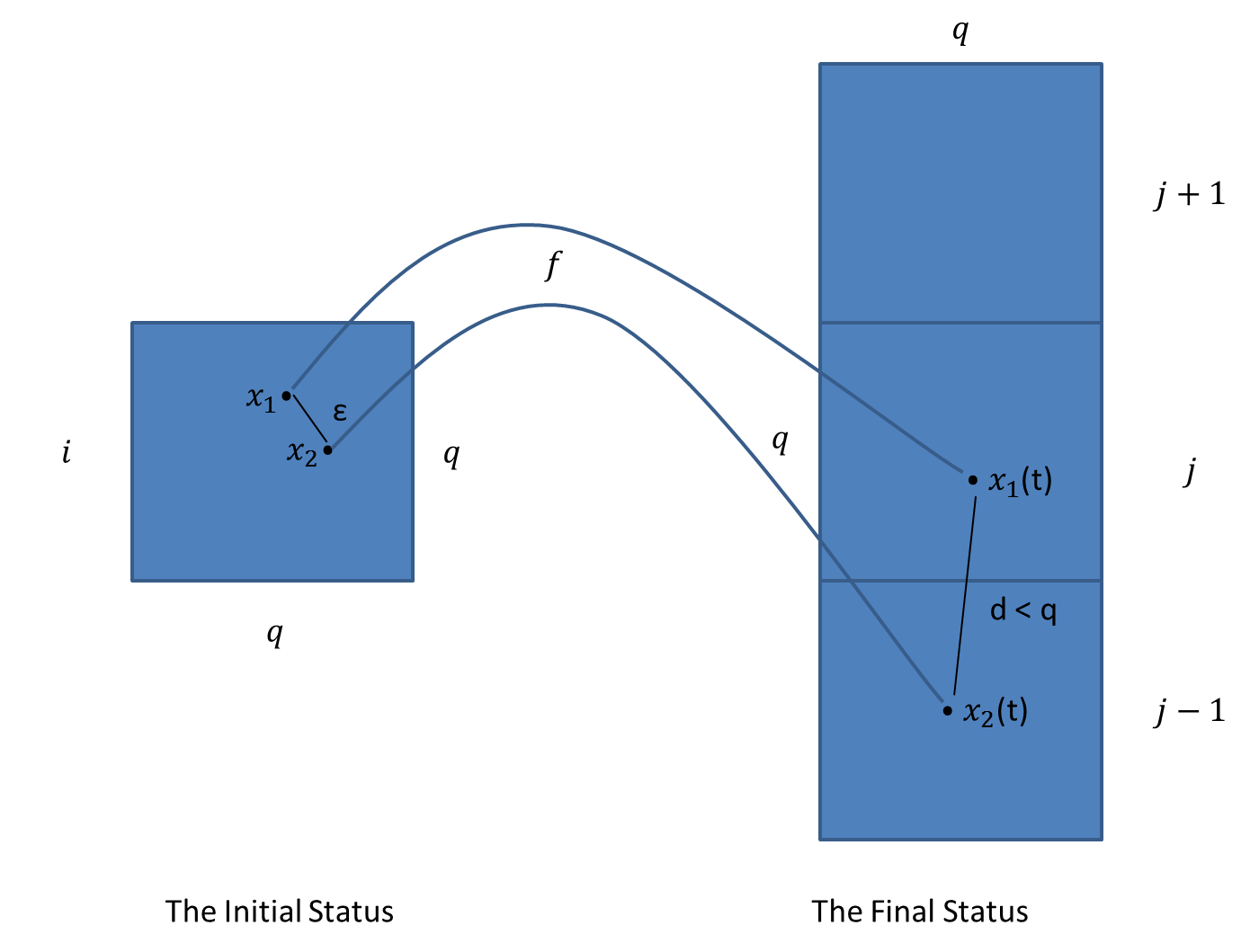}}

  \caption{(a)Two points in box $i$ evolve to boxes $j+1$ and $j-1$ and without further samples, the resulting Ulam-Galerkin matrix would miss a transition $i \to j$ which is expected in a continuous transformation $f$. (b) The success scenario we wish to guarantee. Points $x_1$ and $x_2$ are close enough that the $f$ never casts them further than across successive rectangles $j-1$ and $j$. Recall that a continuous $f$ will map a connected set to a connected set and this should be reflected in the images of the coarse representation of boxes cast across boxes.}
  \label{PA1}
\end{figure}

By Gronwall's inequality,
{\begin{eqnarray}
|x_1(t)-x_2(t)| \le |x_1(t_0)-x_2(t_0)| e^{M|t-t_0|},
\end{eqnarray}}
where $M$ is the Lipschitz constant. Assuming $f$ is uniformly continuously differentable,
{\begin{eqnarray}
M= \max_{ (x,t) \in X \times [t_0,t]} \left | \frac{\partial f}{\partial x} (x,t) \right |.
\end{eqnarray}}
We must control the distance between $x_1(t)$ and $x_2(t)$  through controlling the distance between $x_1(t_0)$ and $x_2(t_0).$ 
Let 
{\begin{eqnarray}
\varepsilon= \frac {q}{e^{M|t-t_0|}}.
\end{eqnarray}}
Then we have,
{\begin{eqnarray}
|x_1(t)-x(t)| \le |x_1(t_0)-x(t_0)| e^{M|t-t_0|} < \varepsilon e^{M|t-t_0|} = q,
\end{eqnarray}}

for all $x \in X$ satisfying $|x_1(t_0)-x(t_0)| < \varepsilon$. That is, for any $x(t_0)$ in an $\varepsilon$-ball of $x_1(t_0)$, $x(t_0)$ will keep its distance with $x_1(t_0)$  less than $q$ through time $t-t_0$.

\item Next, we consider the more general case of  the $\varepsilon$-ball of $x_1(t_0)$; consider $\varepsilon$-balls centered on each sample data $x$ which must cover the whole $X$. 
Since $X$ is compact, the $\varepsilon$-balls of a finite subset $\{ x_k\}_{k=1}^{P}$ cover X.
Let $\{ x_k\}_{k=1}^{P}$ be initial conditions of the flow, with images under the flow $\{ x_k(t)\}_{k=1}^{P}$ in a finite time $t-t_0$.
Define,
$$l= \min_{\{ x_k(t)\}_{k=1}^{P}} \{ \mbox{distance between each of four boundaries of some} \  j \  \mbox{box and} \ x_k(t) \}$$
where the $j$ box is as defined above, is the image of $X$ which we denote $X(t)$. 
See Figure \ref{PA3}.

\begin{figure}[htb!]
  \centering
  \subfloat[some $j$ box from $X(t)$]{\label{fig:gull}\includegraphics[width=.7\textwidth]{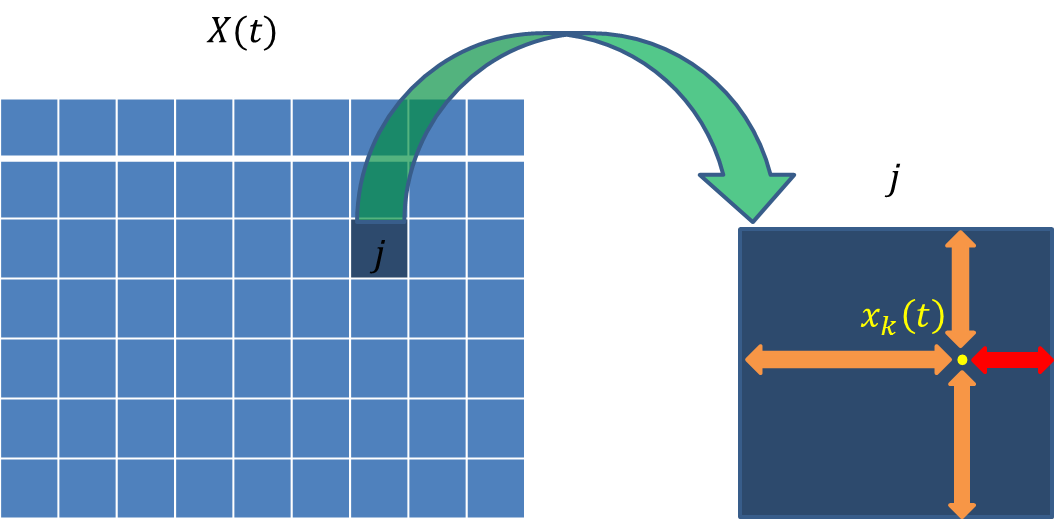}}             
  \caption{The red arrow indicates the minimum distance between a point $x_k(t)$ for some k and four boundaries of a box j for some j that contains $x_k(t)$. We must check all $x_k$ in $\{ x_k(t)\}_{k=1}^{P}$ with all the boxes of the partition of $X(t)$ to develop $l$.}
  \label{PA3}
\end{figure}
Here, again, we use Gronwall's inequality as follows,
$$|x(t)-x_k(t)| \le |x(t_0)-x_k(t_0)| e^{M|t-t_0|} \le l,$$
where $k=1,2, \  ... \ P$.
Let
{\begin{eqnarray}
\varepsilon'= \frac {l}{e^{M|t-t_0|}}, 
\end{eqnarray}}
we have 
{\begin{eqnarray}
|x_k(t)-x(t)| \le |x_k(t_0)-x(t_0)| e^{M|t-t_0|} < \varepsilon e^{M|t-t_0|} = l,
\end{eqnarray}}
which is to say, all points $x$ in an $\varepsilon'$-ball of $x_k$ in $X$ for some $k$ will map into the $l$-ball of $x_k(t)$ in $X(t)$ for the same $k$.

\item In this step, we build a new cover for $X(t)$ with $l$-balls by the following.
There are two types of balls of the cover, Type I covers are those $l$-balls never touch the boundaries of the grid; 
Type II covers are those $l$-balls centered at boundaries and their centers never in the $l$-balls of Type I, see Figure \ref{PA4}.
\begin{figure}[htb!]
  \centering
  \subfloat[two types of covers]{\label{fig:gull}\includegraphics[width=.8\textwidth]{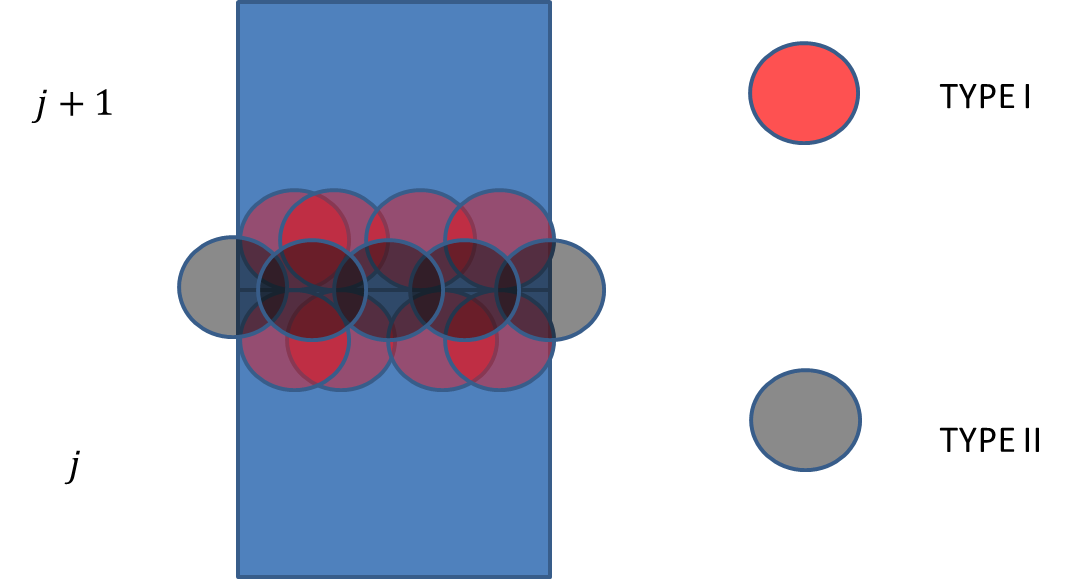}}             
  \caption{}
  \label{PA4}
\end{figure}

$X(t)$ is compact, so we have a finite subset of this new cover still covering $X(t)$. 
Let the centers of these balls be $\{ x_q(t)\}_{q=1}^{Q}$, for Q is some integer.
Since the flow is assumed continuous as is its inverse is continuous, (see \cite{P}), the pre-image of the these $(x_q(t),l)$-balls is still a cover of the pre-image X of $X(t)$, which is also finite.
By the above discussion, more specifically, the finite cover of X consists of $\varepsilon'$-balls of some of $\{ x_q(t)\}_{q=1}^{Q}$. 
That is, the pre-image of Type I $l$-balls in $X(t)$ and another type of covers whose image are Type II $l$-balls in $X(t)$, the shape of these covers may be vary. 
Thus, in $X$, we have the following case, see Figure \ref{PA5}.
\begin{figure}[htb!]
  \centering
  \subfloat[ the pre-image of the two types of covers]{\label{fig:gull}\includegraphics[width=.8\textwidth]{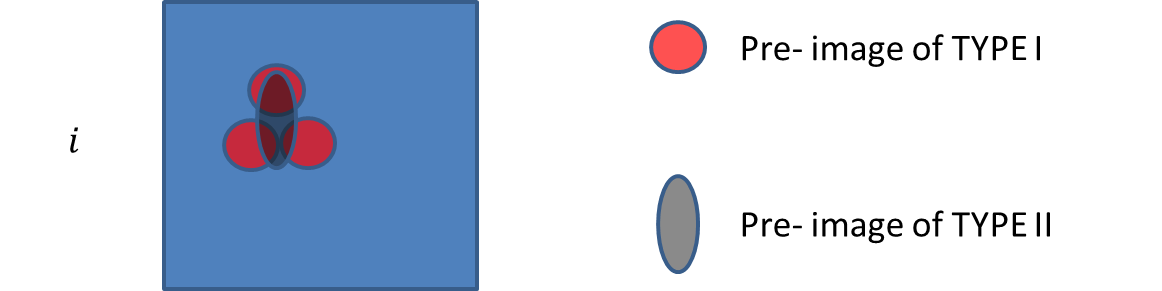}}             
  \caption{We use an ellipse to denote the pre-image of some Type II $l$-ball, but the pre-image may be any shape which depends on the inverse of the flow.}
  \label{PA5}
\end{figure}

Consider the areas which are not overlapping of the balls of both of $X$ and $X(t)$. 
The flow is a bijection between these two areas. 
See Figure \ref{PA6}, which is to say, the points in red $\varepsilon'$-balls in $X$ will never go out of the $l$-balls in $X(t)$.

Notice the important connection between the non-overlapping areas of $X$ and $X(t)$, we can control the non-overlapping area in X by shrink the size of $l$-balls in $X(t)$.
In other words, if $l$ become smaller, so is the non-overlapping area in $X$, that is, the non-overalpping area in $X$ can be ignored for some small $l$.
Then we can guarantee a sufficient sampling for a given grid with length $q$ if we choose mesh-grid of points with distance $\varepsilon'$ or smaller. 
The number of points we need in each box of the partition of the initial status X should be more than $\frac {q}{\varepsilon'}.$ To obtain a $\varepsilon'$, 
we just choose a small $l$ which depends on the Lipschitz constant of a specific problem. 
Note that this description is in terms of a uniform $\varepsilon'$ grid. 
However, a uniformly random cover where this $\varepsilon'$ condition is satisfied is also sufficient.

\begin{figure}[htb!]
  \centering
  \subfloat[ the relationships between the non-overlapping areas]{\label{fig:gull}\includegraphics[width=.6\textwidth]{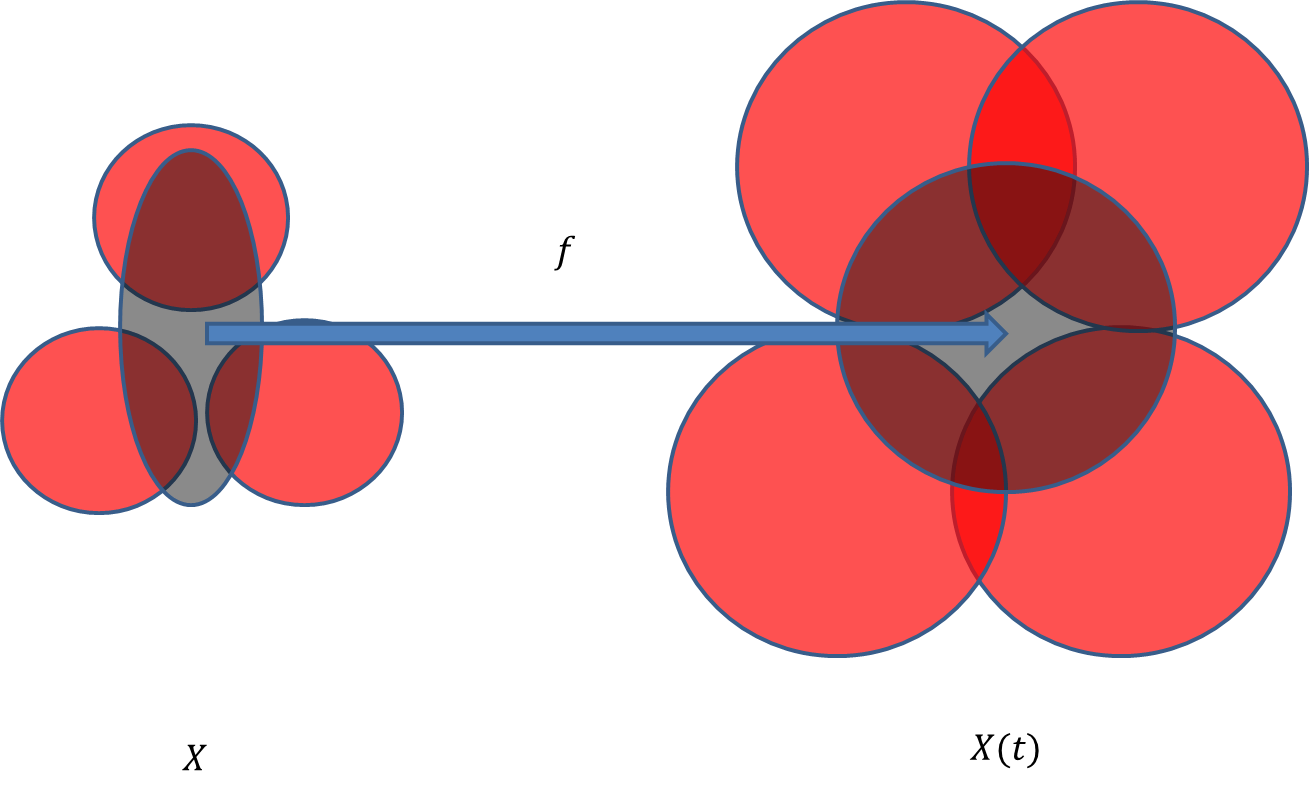}}             
  \caption{}
  \label{PA6}
\end{figure}

\end{enumerate}

%%%%%%%%%%%%%%%%%%%%%%%%%%%%%%%%%%%%%%%%%%%%%%%%%%%%%%%%%%%%%%%%%%%%%%%%%%%%%%%%%
%                                                                                                                    Reference
%%%%%%%%%%%%%%%%%%%%%%%%%%%%%%%%%%%%%%%%%%%%%%%%%%%%%%%%%%%%%%%%%%%%%%%%%%%%%%%%%

%

\end{document}